\documentclass[11pt]{article}
\usepackage{amssymb,amsthm,amsmath}
\usepackage{bbm}

\theoremstyle{definition}
\newtheorem{definition}{Definition}
\theoremstyle{plain}
\newtheorem{lemma}[definition]{Lemma}
\newtheorem{theorem}[definition]{Theorem}

\newcommand{\N}{\mathbb N}
\newcommand{\Q}{\mathbb Q}
\newcommand{\D}{{\mathbb D}} 
\newcommand{\1}{{\mathbbm 1}}
\newcommand\dint{{\rm d}}

\newlength{\fixboxwidth}
\begin{document}
\title{On lower bounds for the $L_2$-discrepancy}
\author{Aicke Hinrichs\thanks{Research of the first author was supported by the DFG Heisenberg grant Hi 584/3-2. } \ and 
        Lev Markhasin\thanks{Research of the second author was supported by a scholarship of the G. C. Starck Foundation.}}
\date{\today}
\maketitle

\begin{abstract}
 The $L_2$-discrepancy measures the irregularity of the distribution of a finite point set.
 In this note we prove lower bounds for the $L_2$ discrepancy of arbitrary $N$-point sets.
 Our main focus is on the two-dimensional case.
 Asymptotic upper and lower estimates of the $L_2$-discrepancy in dimension 2 are well-known and are of the sharp order $\sqrt{\log N}$.
 Nevertheless the gap in the constants between the best known lower and upper bounds is unsatisfactory large for a two-dimensional problem.
 Our lower bound  improves upon this situation considerably.
 The main method is an adaption of the method of K. F. Roth using the Fourier coefficients of the discrepancy function with respect to
 the Haar basis.
 We obtain the same improvement in the quotient of lower and upper bounds in the general $d$-dimensional case.
 Our lower bounds are also valid for the weighted discrepancy.
\end{abstract}

\noindent{\footnotesize {\it 2010 Mathematics Subject Classification.} Primary 11K06,11K38,42C10,65C05.}

\section{Introduction}
 
 The $L_2$-discrepancy is a measure for the irregularity of the distribution of a finite point set with respect to the
 uniform distribution.
 If ${\cal P}$ is an $N$-point subset of the $d$-dimensional unit cube $\Q^d = [0,1)^d$ the {\em discrepancy function} 
 $D_{\cal P}$ is defined as
 \begin{equation}
  \label{eq:_disc}
    D_{\cal P}(x) :=  \sum_{z\in{\cal P}} \1_{C_z}(x) -  \, N \, |B_x|. 
 \end{equation}
 $|B_x|=x_1 \ldots x_d$ denotes the volume of the rectangular box 
 $B_x=[0,x_1)\times \ldots \times [0,x_d)$ for $x=(x_1,\ldots,x_d)\in \Q^d$ and
 $\1_{C_z}$ is the characteristic function of the rectangular box $C_z=(z_1,1)\times\ldots\times(z_d,1)$ for 
 $z=(z_1,\ldots,z_d)\in {\cal P}$. 
 Observe that the sum in this definition is just the number of points of ${\cal P}$ in the box $B_x$. 
 So the discrepancy function measures the deviation of this number from the fair number of points $N \, |B_x|$ which
 would be achieved by a perfect (but impossible) uniform distribution of the points of ${\cal P}$.
 The $L_2$-discrepancy of  ${\cal P}$ is the $L_2$-norm of the discrepancy function given by
 \begin{equation}
  \label{eq:_L2disc}
    \| D_{\cal P} | L_2 \|^2 = \int_{\Q^d} D_{\cal P}(x)^2 \, \dint x . 
 \end{equation}
 
 In this note we are mainly interested in the two-dimensional case $d=2$. Only in the very end we add some
 remarks about the general $d$-dimensional case. So from now on until further notice fix $d=2$. In this case,
 the asymptotic behavior of the minimal possible $L_2$-discrepancy of an $N$-point set for $N\to\infty$ is well-known.
 There is a constant $c$ such that for all $N\in\N$ and all $N$-point subsets ${\cal P} \subset \Q^2$
 \begin{equation}
  \label{eq:_disc_lb}
   \| D_{\cal P} | L_2 \| \ge c  \sqrt{\log N}
 \end{equation} 
 and there is a constant $C$ such that for all $N\in\N$ there exists an  $N$-point subset ${\cal P} \subset \Q^2$ 
 with
 \begin{equation}
  \label{eq:_disc_ub}
   \| D_{\cal P} | L_2 \| \le C  \sqrt{\log N}.
 \end{equation}  
 The lower bound in \eqref{eq:_disc_lb} is the celebrated result by K. F. Roth \cite{R54}.
 Constructions of point sets satisfying \eqref{eq:_disc_ub} are plenty, the first one was given
 by H. Davenport \cite{D56}. 
 For further constructions and the general theory of discrepancy we refer the reader to the books \cite{DP10,M99,NW08,NW10}.
 
 We are interested in the constants in \eqref{eq:_disc_lb} and \eqref{eq:_disc_ub} for large $N$, so let us define
 \begin{equation}
  \label{eq:_defcC}
   \underline{c} = \liminf_{N\to\infty} \inf_{\# {\cal P}=N} \frac{\| D_{\cal P} | L_2 \|}{\sqrt{\log N}}   
   \ \ \ \mbox{and}\ \ \ 
   \overline{c} = \limsup_{N\to\infty} \inf_{\# {\cal P}=N} \frac{\| D_{\cal P} | L_2 \|}{\sqrt{\log N}}   
 \end{equation}  
 The best estimates for the constants $\underline{c}$ and $\overline{c}$ known so far are
 \begin{equation}
  \label{eq:_bestc}
   0.0046918\ldots = \sqrt{\frac{1}{2^{16}\log 2}} \ \le \  \underline{c} \ \le \ \overline{c} \ \le \sqrt{\frac{278629}{2811072 \, \log 22}} =  0.17907\ldots
 \end{equation}  
 The bound for $\underline{c}$ is from a modification of the proof of Roth and can be found in \cite[Chapter 2, proof of Lemma 2.5]{KN74}.
 With this constant the estimate \eqref{eq:_disc_lb} even holds for all $N\in\N$.
 The bound for $\overline{c}$ is from a recent construction in \cite{FPPS10} using generalized scrambled Hammersley point sets.
 For a two-dimensional problem, the gap between the constants is huge. 
 
 The main purpose of this note is to improve the lower bound. Our main result is
 \begin{theorem}
  \label{thm:_main}
  For all $N\in\N$ and all $N$-point sets ${\cal P} \subset \Q^2$ the inequality
  $$ \| D_{\cal P} | L_2 \| \ge c  \sqrt{\log N} $$
  holds with
  $$ c = {\frac{7}{216 \, \sqrt{\log 2}}} = 0.038925\ldots.$$
 \end{theorem}

 The proof is still a variant of the method of Roth which uses the information that certain dyadic rectangles do not
 contain any points of a given point set and adds those local discrepancies up with the help of orthogonal functions.
 Our improvement is due to the fact that we consider different levels of dyadic rectangles. A convenient method to do
 this is to compute Fourier coefficients with respect to the Haar system and then use Parseval's formula. This method
 was already used in a recent paper \cite{H10} of the first author to prove optimal upper estimates for the discrepancy of
 Hammersley type point sets measured in spaces of dominating mixed smoothness.
 
 The next section contains the necessary tools concerning the Haar basis in $L_2$. In Section 3 we prove the
 lower bound. In the final section we show that our lower bound remains valid for the weighted discrepancy
 and we comment on what can be done with the Haar function method in higher dimensions.
 
\section{Haar coefficients of the discrepancy function}

 A dyadic interval of length $2^{-j}, j\in \N_0,$ in $[0,1)$ is an interval of the form 
 $I=I_{j,m}:=\big[ 2^{-j}m,2^{-j}(m+1)\big)$ for $m=0,1,\ldots,2^j-1$. 
 The left and right half of $I=I_{j,m}$ are the dyadic intervals $I^+ = I_{j,m}^+ =I_{j+1,2m}$ and
 $I^- = I_{j,m}^- =I_{j+1,2m+1}$, respectively. 
 The Haar function $h_I = h_{j,m}$ with support $I$ is the function on $[0,1)$ which is 
 $+1$ on the left half of $I$, $-1$ on the right half of $I$ and 0 outside of $I$. 
 The $L_\infty$-normalized Haar system consists of all Haar functions $h_{j,m}$ with $j\in\N_0$ and  
 $m=0,1,\ldots,2^j-1$ together with the indicator function $h_{-1,0}$ of $[0,1)$.
 Normalized in $L_2(\Q)$ we obtain the orthonormal Haar basis of $L_2(\Q)$. 

 Let $\N_{-1}=\{-1,0,1,2,\ldots\}$ and define $\D_j=\{0,1,\ldots,2^j-1\}$ for $j\in\N_0$ and $\D_{-1}=\{0\}$ for $j=-1$.
 For $j=(j_1,j_2)\in\N_{-1}^2$ and $m=(m_1,m_2)\in \D_j :=\D_{j_1}\times \D_{j_2}$, the Haar function $h_{j,m}$
 is given as the tensor product $h_{j,m} (x) = h_{j_1,m_1}(x_1)\, h_{j_2,m_2}(x_2)$ for $x=(x_1,x_2)\in[0,1)^2$.
 We will call the rectangles $I_{j,m} = I_{j_1,m_1} \times  I_{j_2,m_2}$ dyadic rectangles.
 The $L_\infty$-normalized tensor Haar system consists of all Haar functions $h_{j,m}$ with $j\in\N_{-1}^2$ and  
 $m \in \D_j$. Normalized in $L_2(\Q^2)$ we obtain the orthonormal Haar basis of $L_2(\Q^2)$. 
 
 Now Parseval's equation shows that the $L_2$-norm of a function $f\in L_2(\Q^2)$ can be computed as
 \begin{equation}
  \label{eq:_Parseval}
   \|f|L_2\|^2 =  \sum_{j\in\N_{-1}^2}  2^{\max(0,j_1)+\max(0,j_2)}  \sum_{m\in \D_j} |\mu_{j,m}|^2  
 \end{equation}  
 where 
 \begin{equation}
  \label{eq:_haarcoeffs}
   \mu_{j,m} = \mu_{j,m}(f) = \int_{\Q^2} f(x) h_{j,m}(x) \, \dint x
 \end{equation}
 are the Haar coefficients of $f$.
 
 The following two crucial lemmas are easy to verify and were already used in \cite{H10}.
\begin{lemma}
 \label{lem:_1}
 Let $f(x)=x_1\,x_2$ for $x=(x_1,x_2)\in \Q^2$. Let $j\in\N_0^2$, $m\in \D_j$ and let $\mu_{j,m}$ be the 
 Haar coefficient of $f$ given by \eqref{eq:_haarcoeffs}.
 Then 
 $$\mu_{j,m}=2^{-2 j_1 -2 j_2 - 4}.$$
\end{lemma} 

\begin{lemma}
 \label{lem:_2}
 Fix $z=(z_1,z_2)\in\Q^2$ and let $f(x)=\1_{C_z}(x)$ for $x=(x_1,x_2)\in \Q^2$. Let $j\in\N_0^2$, $m\in \D_j$ and let $\mu_{j,m}$ be the 
 Haar coefficient of $f$ given by \eqref{eq:_haarcoeffs}. Then $\mu_{j,m}=0$ whenever $z$ is not contained in the interior of the dyadic rectangle
 $I_{j,m}$ supporting $h_{j,m}$.
\end{lemma}

\section{The lower bound}

 We are now ready to prove Theorem \ref{thm:_main}. 
 Let $N\in\N$ with $N\ge 2$ and let  ${\cal P} \subset \Q^2$ be an $N$-point set.
 Let $j=(j_1,j_2)\in \N_{0}^2$, $m\in \D_j$ be such that no point of  ${\cal P}$ lies in the interior of the dyadic rectangle
 $I_{j,m}$ supporting $h_{j,m}$.
 Let $\mu_{j,m}$ be the Haar coefficient of the discrepancy function \eqref{eq:_disc}.
 Now Lemmas \ref{lem:_1} and \ref{lem:_2} imply that 
 $$ \mu_{j,m} = - N \, 2^{-2 j_1 -2 j_2 -4}. $$

 Observe that for fixed $j=(j_1,j_2)\in \N_{0}^2$ the cardinality of $\D_j$ is $2^{j_1+j_2}$ and the interiors of the dyadic rectangles  
 $I_{j,m}$ supporting $h_{j,m}$ are mutually disjoint. This implies that there are at least $2^{j_1+j_2}-N$ such $m\in \D_j$ for which
 no point of  ${\cal P}$ lies in the interior of the dyadic rectangle $I_{j,m}$ supporting $h_{j,m}$.
 
 We abbreviate $M=\lceil \log_2 N \rceil$. Then we obtain from \eqref{eq:_Parseval} that
 \begin{eqnarray*}
   \| D_{\cal P} | L_2 \|^2 
   &\ge& 
   N^2 \sum_{j_1+j_2 \ge M} 2^{j_1+j_2} (2^{j_1+j_2}-N)  2^{-4 j_1 -4 j_2 - 8}\\
   &=& 2^{-8} N^2  \sum_{j_1+j_2 \ge M} 4^{-(j_1+j_2)} - 2^{-8} N^3 \sum_{j_1+j_2 \ge M} 8^{-(j_1+j_2)}. \\
 \end{eqnarray*}
 Now for any $q>1$ we have
 $$ \sum_{j_1+j_2 \ge M} q^{-(j_1+j_2)} = \sum_{k=M}^\infty (k+1) q^{-k} = q^{-M+1} \Big( \frac{M}{q-1} + \frac{q}{(q-1)^2} \Big)$$
 which leads to
 \begin{eqnarray*}
   \| D_{\cal P} | L_2 \|^2 
   &\ge& 
   2^{-6} (N 2^{-M})^2  \Big( \frac{M}{3} + \frac{4}{9} \Big)
   -
   2^{-5} (N 2^{-M})^3  \Big( \frac{M}{7} + \frac{8}{49} \Big) \\
   &\ge& 
   2^{-6} (N 2^{-M})^2  \,  \frac{M}{3} 
   -
   2^{-5} (N 2^{-M})^3  \,  \frac{M}{7}  \\
 \end{eqnarray*} 
 where the last estimate easily follows from $0< N 2^{-M} \le 1$.
 
 Let now $t=M - \log_2 N$ so that $0\le t< 1$ and $N 2^{-M} = 2^{-t}$. Then we have proved 
 $$ \| D_{\cal P} | L_2 \|^2  \ge \gamma \, \log_2 N $$ 
 if we can verify that
 $$ 2^{-6} 2^{-2t}  \,  \frac{M}{3} - 2^{-5} 2^{-3t}  \,  \frac{M}{7} \ge \gamma (M-t) $$ 
 for all $M\in \N$ and $0\le t< 1$. The last inequality is equivalent to
 $$ \left( \gamma - 2^{-6}\, 3^{-1} 2^{-2t} +  2^{-5}\, 7^{-1} 2^{-3t} \right)\, M \le \gamma \, t $$
 which is certainly satisfied whenever $\gamma\ge 0$ and 
 $$ \gamma \le  2^{-6}\, 3^{-1} 2^{-2t} -  2^{-5}\, 7^{-1} 2^{-3t} $$
 for all $0\le t < 1$ or, alternatively, 
 $$ \gamma \le  2^{-6}\, 3^{-1} y^2 -  2^{-5}\, 7^{-1} y^3 $$
 for all $1/2< y\le 1$. The maximal value of the right hand side is easily seen to be $\frac{49}{46656}$ for $y=\frac{7}{9}$.
 So we can choose 
 $$ \gamma = \frac{49}{46656} = \Big( \frac{7}{216} \Big)^2 $$
 which leads to the value for $c$ in the theorem. This finishes the proof. \qed
 
\section{Final Remarks}

Our lower bound is also valid for the weighted discrepancy which can be defined as follows. 
Let $a=(a_z)_{z\in\cal P}$ be a system of real numbers associating a weight $a_z$ with a point $z\in\cal P$.
Then the weighted discrepancy function is defined as
$$    D_{{\cal P}, a}(x) :=  \sum_{z\in{\cal P}} a_z\1_{C_z}(x) -  \, N \, |B_x|.$$
The discrepancy function defined by \eqref{eq:_disc} is obtained in the case that all points of $\cal P$ have weight $1$. 
Thanks to Lemma \ref{lem:_2} the Haar coefficient with respect to a Haar function whose support does not intersect $\cal P$
does not depend on the weights.
So one gets the same lower bound with the same constant for the weighted $L_2$-discrepancy as in the case without weights.
Hence we have the following generalization of Theorem \ref{thm:_main} to the weighted discrepancy.
 \begin{theorem}
  \label{thm:_main2}
  For all $N\in\N$, all $N$-point sets ${\cal P} \subset \Q^2$ and all weights $a=(a_z)_{z\in\cal P}$ the inequality
  $$ \| D_{{\cal P},a} | L_2 \| \ge c  \sqrt{\log N} $$
  holds with
  $$ c = {\frac{7}{216 \, \sqrt{\log 2}}} = 0.038925\ldots.$$
 \end{theorem}

We now consider point sets in higher dimensions. Then for all $N\in\N$ and all $N$-point subsets ${\cal P}\subset\Q^d$ there is a known lower bound of the form
\[\| D_{\cal P} | L_2 \| \ge c_d  (\log N)^\frac{d-1}{2}\]
where the constant is known from \cite{DP10} as
\[c_d=\frac{1}{2^{2d+4}\sqrt{(d-1)!}(\log 2)^\frac{d-1}{2}}.\]

Our intention is to use the method from above to improve this constant. The idea of tensor product Haar bases can be easily transfered to higher dimensions $d> 2$. For $j=(j_1,\ldots,j_d)\in\N_{-1}^d$ and $m=(m_1,\ldots,m_d)\in \D_j:=\D_{j_1}\times\ldots\times \D_{j_d}$ the Haar function $h_{j,m}$ is given as the tensor product $h_{j,m}(x)=h_{j_1,m_1}(x_1)\ldots h_{j_d,m_d}(x_d)$ for $x=(x_1,\ldots,x_d)\in [0,1)^d$. We will call the rectangles $I_{j,m}=I_{j_1,m_1}\times\ldots\times I_{j_d,m_d}$ dyadic boxes. The $L_\infty$-normalized tensor Haar system consists of all Haar functions $h_{j,m}$ with $j\in\N_{-1}^d$ and  
 $m \in \D_j$. Normalized in $L_2(\Q^d)$ we obtain the orthonormal Haar basis of $L_2(\Q^d)$.

Now Parseval's equation shows that the $L_2$-norm of a function $f\in L_2(\Q^d)$ can be computed as
 \begin{equation}
  \label{eq:_Parseval_d}
   \|f|L_2\|^2 =  \sum_{j\in\N_{-1}^d}  2^{\max(0,j_1)+\ldots+\max(0,j_d)}  \sum_{m\in \D_j} |\mu_{j,m}|^2  
 \end{equation}  
 where
 \begin{equation}
  \label{eq:_haarcoeffs_d}
   \mu_{j,m} = \mu_{j,m}(f) = \int_{\Q^d} f(x) h_{j,m}(x) \, \dint x.
 \end{equation}
 are the Haar coefficients of $f$.
 
Analogously to the case $d=2$ we can state the following two lemmas. They are easy to verify.

\begin{lemma}
 \label{lem:_3}
 Let $f(x)=x_1\ldots x_d$ for $x=(x_1,\ldots,x_d)\in \Q^d$. Let $j\in\N_0^d$, $m\in \D_j$ and let $\mu_{j,m}$ be the 
 Haar coefficient of $f$ given by \eqref{eq:_haarcoeffs_d}.
 Then 
 $$\mu_{j,m}=2^{-2 j_1-\ldots -2 j_d - 2d}.$$
\end{lemma} 

\begin{lemma}
 \label{lem:_4}
 Fix $z=(z_1,\ldots,z_d)\in\Q^d$ and let $f(x)=\1_{C_z}(x)$ for $x=(x_1,\ldots,x_d)\in \Q^d$. Let $j\in\N_0^2$, $m\in \D_j$ and let $\mu_{j,m}$ be the 
 Haar coefficient of $f$ given by \eqref{eq:_haarcoeffs_d}. Then $\mu_{j,m}=0$ whenever $z$ is not contained in the interior of the dyadic box
 $I_{j,m}$ supporting $h_{j,m}$.
\end{lemma}

Now let $N\in\N$ with $N\ge 2$ and let  ${\cal P} \subset \Q^d$ be an $N$-point set.
 Let $j\in\N_0^d$, $m\in \D_j$ be such that no point of  ${\cal P}$ lies in the interior of the dyadic box
 $I_{j,m}$ supporting $h_{j,m}$.
 Let $\mu_{j,m}$ be the Haar coefficient of the discrepancy function \eqref{eq:_disc}.
 Then Lemmas \ref{lem:_1} and \ref{lem:_2} imply that 
 $$ |\mu_{j,m}| = N \, 2^{-2 j_1-\ldots -2 j_d -2d} $$
 for $j=(j_1,\ldots,j_d)\in \N_{0}^2$. 
 We obtain the following result.
 
 \begin{theorem}
  \label{thm:_main_d}
  For all $N\in\N$ and all $N$-point sets ${\cal P} \subset \Q^d$ the inequality
  $$ \| D_{\cal P} | L_2 \| \ge \;c_d  (\log N)^\frac{d-1}{2} $$
  holds with
  $$ c_d = {\frac{7}{27\cdot 2^{2d-1}\,\sqrt{(d-1)!} \, (\log 2)^{\frac{d-1}{2}}}}.$$
 \end{theorem}
 
In comparison with the known result the constant is improved by a factor of $\frac{224}{27}=8.296296\ldots$. The calculation of the constant is analogous to the case $d=2$. First one obtains
\[\| D_{\cal P} | L_2 \|^2\geq 2^{-4d} N^2  \sum_{j_1+\ldots+j_d \ge M} 4^{-(j_1+\ldots+j_d)} - 2^{-4d} N^3 \sum_{j_1+\ldots+j_d \ge M} 8^{-(j_1+\ldots+j_d)}.\]
Then one checks that the coefficient of $M^{d-1}$ in
\[\sum_{j_1+\ldots+j_d \ge M} q^{-(j_1+\ldots+j_d)}\]
for any $q>1$ is 
\[\frac{q^{-M+1}}{(q-1)(d-1)!}.\]
Finally one obtains
\[\| D_{\cal P} | L_2 \|^2 \geq 2^{-4d} (N 2^{-M})^2\frac{4}{3}\frac{M^{d-1}}{(d-1)!}-2^{-4d} (N 2^{-M})^3\frac{8}{7}\frac{M^{d-1}}{(d-1)!}.\]
Then analogously to the 2-dimensional case we get the estimate
$$ \| D_{\cal P} | L_2 \|^2  \ge \gamma \, (\log_2 N)^{d-1} $$ 
if
$$ \gamma \le  \frac{1}{2^{4d}(d-1)!}\left(\frac{4}{3} y^2 -  \frac{8}{7} y^3\right) $$
for all $1/2< y\le 1$. The maximal value of the right hand side is reached for $y=\frac{7}{9}$ and is
\[\gamma=\frac{1}{2^{4d}(d-1)!}\left(\frac{14}{27}\right)^2.\]
This leads to the value of $c_d$ in the theorem.

Analogously to the $2$-dimensional case one obtains in the $d$-dimensional case the same bounds for the weighted discrepancy as for the unweighted.

\vspace{10mm}

\noindent
Mathematisches Institut, Friedrich-Schiller-Universit\"at Jena,\\ 
Ernst-Abbe-Platz 2, D-07737 Jena, Germany\\
E-mail address:  a.hinrichs@uni-jena.de, lev.markhasin@uni-jena.de
\end{document}